\documentstyle[a4,12 pt,amssymb]{article}

\newcommand{\vv}{\vspace{2mm}}
\newcommand{\Gal}{\mbox{{\rm Gal}}}
\newcommand{\fp}{{\frak p}}
\def\C{{\cal C}}
\def\Q{{\bf Q}}

\def\F{{\bf F}}
\def\P{{\bf P}}

\def\ext{\wedge^r}
\newcommand{\YY}{\tilde{Y}}

\def\Spec{{\mathrm Spec}}
\def\pf{{\indent\textit{Proof.}\ }}
\def\qed{\hfill$\square$}
\newcounter{para}[section]
\renewcommand{\thepara}{\thesection.\arabic{para}}
\renewcommand{\thesection}{\arabic{section}}
\renewcommand{\paragraph}{\refstepcounter{para}
\indent{\bf{\thepara}}}
\newcommand{\sectioning}{\refstepcounter{section}
\indent{\bf \thesection.}}

\newcommand{\n}{{\mathfrak n}}
\renewcommand{\a}{{\mathfrak a}}

\begin{document}

{\footnotesize \noindent Running head: Drinfeld modules and Galois groups \hfill September 3, 2002

\noindent Math.\ Subj.\ Class.\ (2000): 12F12, 11G09}

\vspace{45mm}

\begin{center} 

{\Large Torsion of Drinfeld modules and

\vv

equicharacteristic unimodular Galois covers}

\vv

{\sl by} Gunther Cornelissen {\sl and} Marina Tripolitaki

\vv

\end{center}

\vv

\vv

{\bf Introduction} 

\vv

\vv

\noindent A basic problem of inverse Galois theory is to 
find, for a given field $K$ and finite group $G$, an absolutely irreducible cover of the
projective line ${\bf P}^1_K$ defined over $K$, whose Galois group 
equals $G$. The outstanding case is $K=\Q$ (viz., solving the inverse Galois theory problem over ${\bf Q}(t)$), but here we will be concerned
with its cousin $K={\bf F}_q(T)$, a global function field over a finite field $\F_q$
with $q$ elements. A famous theorem of Harbater-Pop (\cite{Pop}) says that any $G$ can be realized over ``large'' fields (e.g., the separable closure of $K$), but the question remains over what smaller field such a cover can be defined (especially in connection with the so-called ``arithmetic Abhyankar conjecture'', cf.\ \cite{Abhyankar:95}, (9.2)). 

In this direction, we have the following result of Nori (\cite{Nori}): if $G=A(\F_{q^d})$ is the group of $\F_{q^d}$ points of a semi-simple simply connected algebraic 
group $A$ 
over $\F_q$, then $G$ is a Galois group 
over $K$ (ramified above at most one point) --- this applies in particular to $SL(r,q)$, so there is also such a $PSL(r,q^d)$-cover. However, this gives only one such cover, not a ``family'' over ${\bf P}^1_K$ for $K$ a {\sl global} field. In this paper, we link this problem with torsion of Drinfeld modules to show that there
is at least an $(r-1)$-dimensional family of $K$-regular unimodular covers (which can
then be given {\sl explicitly} by calculations with modular forms):

\vv

{\bf Theorem.} \ \ {\sl Let $q$ be a power of an odd prime number $p$. For every pair of integers $(r,d)$, such that $r$ is even and $(r;q^d-1)=2$ there exists a $K$-regular Galois extension of the purely transcendental field $K(t_1,...,t_{r-1})$ with Galois group the unimodular group $PSL(r,q^d)$.}

\vv

To appreciate this result, let us describe two known related facts. The general linear group $GL(m,q^d)$ (and, {\sl a fortiori}, its projectivized version) is known to occur as a Galois group over $K(t_1,...,t_{r-1})$, since it occurs as the image of the Galois representation $\rho_{\phi,\fp,K}$ associated
to the ${\frak p}$-torsion of the ``generic'' rank-$r$ Drinfeld module $\phi$ over $K$, where ${\frak p}$ is irreducible 
of degree $d$ in ${\bf F}_q[T]$ (but note that these extensions are in general not regular, since they contain the $\frak p$-th ray class field of $K$, cf.\ infra); explicit
proofs are due to Abhyankar and collaborators (cf.\ \cite{Abhyankar:00b}). On the other hand, the unimodular group $PSL(2,p)$ is known to be the Galois group
of an extension of ${\bf Q}(t)$ if $2,3$ or $7$ is not a square modulo $p$, by the work of K.-y. Shih (\cite{Shih:78}). In modern language, he uses moduli of elliptic ${\bf Q}$-curves. 
 
As a matter of fact, our proof of the above theorem will be a combination
of the techniques used to prove these two results: we study a twisted version of moduli spaces
of Drinfeld modules of rank $r$ with an appropriate level structure. The main novelties are:
(a) the fact that we use higher dimensional varieties -- hence rationality questions related to moduli spaces become more difficult; (b) the fact that the Weil pairing is replaced by determinants in the category of $t$-motives; (c) the fact that there is no residue class condition in the final result.

Here is a more detailed outline: we start by looking at  the moduli scheme $Y^{r}_0(T)$
 of Drinfeld modules over $K$ of rank $r$ ``with level $\Gamma_0(T)$-structure'' as classifying Drinfeld modules of rank $r$ with a ``full flag'' of subgroups of its $T$-torsion. We give an explicit description of this variety.  Let $L=K(\sqrt{\fp})$ for $\fp$ prime in $\F_q[T]$ of degree $d$.
We define $\YY$ to be the quotient $$\YY := (\langle w \rangle \times \Gal(L/K) )\backslash (Y^{r}_0(T) \times_{\Spec\ K} \Spec\ L),$$ where $w$ is an Atkin-Lehner-style involution (whose action we describe explicitly) -- this will turn
out to make sense since $r$ is even. We then show that $\YY$  has a rational function field by showing that a naive compactification is a Brauer-Severi variety with 
a $K$-rational point. If $(q^{\deg \fp}-1,r)=2$, then to any $K$-rational point of $\YY$, there is associated a representation 
$$ \rho : \Gal(\bar{K}/K) \rightarrow PGL(r,{\bf F}_q[T]/\fp), $$
whose determinant is minus the quadratic character of  $\zeta T$ modulo $\fp$ for a constant $\zeta \in \F_q^*$ that only depends on the rank $r$. We then show that this character can take any value for a suitable choice of $N$:

\vv

{\bf Theorem} (continuation). \ {\sl There exists a computable constant 
$\zeta \in \F_q^*$ only depending on $r$, such that if $\fp$ is irreducible
in $\F_q[T]$ of degree $d$ with $\fp(0)\zeta$ a non-square in ${\bf F}_q^*$, then the desired cover is the Galois closure of $Y^r_0(\fp) \times_{Y^r(1)} \YY \rightarrow \YY$; furthermore, such $\fp$ exist.}

\vv 

As this sketch of the proof indicates, the construction can be made explicit if one can find explicit equations for suitable moduli schemes. For example if $r=2$, one has enough control over modular forms (and hence, function fields of modular curves), and one arrives at the following result. 
\vv

{\bf Proposition.} \ {\sl Let $h$ be a modular form of weight $q+1$ and character $\det$ for $GL(2,\F_q[T])$ (which is then unique up to a constant). Let $N \in \F_q[T]$ be
a monic irreducible polynomial of degree two whose constant term is a square. There exists a polynomial $P \in K[x,y]$ of the form $$P(x,y) = T^{q+1} y^2 x^2 + a(y) x + y$$ of bidegree $(2,q^2+1)$ such that $$P(h(Tz)/h(z),\sqrt[q^2-1]{f_T(Nz)/f_T(z)}) = 0.$$ Then the splitting field of the numerator of
$$P(T^{-\frac{q+1}{2}} \frac{\sqrt{N}+x}{\sqrt{N}-x},\frac{\sqrt{N}+y}{\sqrt{N}-y})$$ defines a $K$-regular extension of $K(x)$ with Galois group $PSL(2,q^2)$.}

\vv

The proposition is entirely algorithmical since one can determine the
polynomial $P$ by looking at series expansions of modular forms.

\vv

{\bf Example.} The polynomial in ${\bf F}_3(T)(x)[y]$ given by
{\footnotesize $$ 
T^2x^2y^{10}-T^2Nxy^9-N^2(x^2-N)(Ty^8-Ny^6-N^2y^4-TN^3y^2)+T^2N^5xy-T^2N^6
$$} with $N=T^2+1$ defines a regular extension of ${\bf F}_3(T)(x)$ with Galois group $PSL(2,9)$. 

\vv

{\bf Remark.} The branching of the cover in the proposition is above
the {\sl two} cusps of $\YY$, whereas Nori's covers are unramified 
over ${\bf A}^1$. 

\vv

\vv

\sectioning\label{DM}\ 
{\bf Drinfeld modular schemes}

\vv

\vv

In this section we review the main definitions and properties from 
the theory of moduli spaces of Drinfeld modules, cf.\ 
\cite{Vladuts:91}. Let us note from the start that we need only naive
\'etale level structures, as we don't bother about bad fibers of
the moduli spaces we will consider. 

\vspace{1ex}
\paragraph\label{DM-not} {\bf Notations.} \  Let $\F_q$ be a finite field with $q$ elements ($q$ a power of an odd prime), $A=\F_q[T]$
the polynomial ring over $\F_q$, $K=\F_q(T)$ the field of rational functions,
$K_\infty=\F_q((T^{-1}))$ its completion for the valuation $|a|=q^{\deg(a)}$ for
$a \in A$ and $C$ the completion of an algebraic closure of $K_\infty$.

\vspace{1ex}
\paragraph\label{DM-DM} {\bf Drinfeld modules.} \ For any field $L$ containing $A$, a {\sl Drinfeld module} $\phi$ over $L$ is a non-trivial representation $\phi$ of $A$ 
in the endomorphism ring of the additive group scheme ${\bf G}_{a/L}$ that
induces the identity on the tangent space at zero.  Said otherwise, 
it is a ring morphism $$\phi \ : \ A \rightarrow L\{ \tau \} \ : \ a \rightarrow \phi_a$$ from $A$ to 
the twisted polynomial ring $L\{\tau \}$, where $\tau$ is the $q$-Frobenius
element, {\sl i.e.}, satisfies $\tau x = x^q \tau$ for all $x \in L$. Furthermore, $\phi$ 
satisfies $\phi_a(0)=a$ and $\phi(A) \not \subset L$. It
turns out that $\deg_\tau (\phi_a) = \deg(a) \cdot r$ for some constant
integer $r$, which is called the {\sl rank} of $\phi$. 
 A {\sl morphism} between two Drinfeld modules $\phi$ and $\psi$ over $L$ is 
an element $u \in \mbox{End}({\bf G}_{a/L})$ such that $u\phi=\psi u$. 

A Drinfeld module $\phi$ induces a new $A$-module structure on ${\bf G}_{a/L}$ given 
by $a \cdot x = \phi_a(x)$ for $a \in A, x \in L$, so it makes sense
to speak of the {\sl $\n$-torsion} $\phi[\n]$ (as a subgroup scheme of ${\bf G}_{a/L}$) of $\phi$ for any $\n \in A$. It turns out that $\phi[\n] \cong 
(\n^{-1}A/A)^r$. The $\n$-torsion comes naturally with an action of the Galois group of the separable closure $\bar{K}$ of $K$, so that we get an (equicharacteristic) Galois representation
$$ \rho_{\phi,\n,K} : G_K \rightarrow \mbox{Aut} \ \phi[\n] = GL(r,A/\n), $$
where $G_K$ is the absolute (separable) Galois group of $K$.

\vspace{1ex}
\paragraph\label{DM-moduli} {\bf Moduli spaces with \'etale level structures.} \ It makes sense to speak of 
a Drinfeld module over an $A$-scheme $S$: it consists of a pair consisting
of a line bundle $\mathcal L$ on $S$ and a map $\phi$ from $A$ to
the endomorphisms of $\mathcal L$ as an $S$-group scheme such that the differential 
of $\phi(a)$ equals $a$ for all $a \in A$ and such that for any field $k$
with a morphism $\mbox{Spec}(k) \rightarrow S$, the base change of $\phi$ to $k$ is an ordinary 
Drinfeld module. Let $\n \in A$. An {\sl \'etale level $\n$ structure} on
$(\mathcal{L},\phi)$ is an isomorphism of group schemes $$ \lambda : (\n^{-1} A /A)^r \times_S T \rightarrow \phi[\n] \times_S T$$ (for some finite \'etale base change $T \rightarrow S$) commuting with the $A$-action.
The functor $\mbox{\sf{Sch}}/A \rightarrow \mbox{\sf{ Sets}}$ 
that associates to an $A$-scheme $S$ the set of $S$-Drinfeld modules of rank $r$ with
\'etale level $\n$-structure up to isomorphism is representable by a scheme 
$M^r(\n)$ of finite type over $A$ (here, isomorphism is given locally by an isomorphism of Drinfeld modules compatible with level structures). Seen over $C$, $M^r(\n)$ is the disjoint
union of irreducible components parametrized by $ (A/\n)^*/\F_q^*$, and these
irreducible components are isomorphic over $C$ to the rigid analytic
space $\Gamma(\n) \backslash \Omega^r$, where $\Omega^r$ is Drinfeld's 
$r$-dimensional rigid analytic space and $\Gamma(\n)$ is the principal
congruence group defined by $\Gamma(\n) =\{ \gamma \in GL(r,A) : \gamma = {\bf 1} \mbox{ mod } \n \}.$

\vspace{1ex}
\paragraph\label{CFT} {\bf Rank one Drinfeld modules and explicit class field theory.} \ For $r=1$,
the scheme $M^1(\n)$ is the spectrum of the integral closure 
of $A[\n^{-1}]$ in the $\n$-th ray class field $K_+(\n)$ of $K$ (i.e., $K_+(\n)$ is the fixed
field  of the maximal abelian extension of $K$ that is totally split
at $\infty$ by the $1$-units modulo $\n$). The
natural action of $(A/\n)^*/\F_q^*$ on $M^1(\n)$ is given by class field theory. 
This means the following: let $K(\n)$ be the field extension of $K$ given
by adjoining the $\n$-torsion of any rank one Drinfeld module $\psi$ over $K$ to $K$; then $K_+(\n)$ is its maximal ``real'' subfield (this does not depend
on the choice of $\psi$). Then there is an explicit isomorphism 
\begin{eqnarray*} (A/\n)^* &\rightarrow& \Gal(K(\n)/K) \\ {\frak a} &\mapsto& [ \alpha \in \psi[\n] \mapsto \psi_{\frak a}(\alpha)] \end{eqnarray*}
such that $\psi_{\frak a} (\alpha) = \alpha^{q^{\deg({\frak a})}} \mbox{ mod } \n $ (i.e., ${\frak a} \mapsto \mbox{Frob}_{\frak a}$).

As before, the scheme $M^1(\n)$ splits over $K_+(\n)$ in a number of absolutely irreducible copies of ${\Spec} \ K_+(\n)$ parametrized by $(A/\n)^*/\F_q^*$. 
Note that if we denote by $Y^r(\n)$ a fixed irreducible component of $M^r(\n)$ over
$K_+(\n)$, $Y^r(\n)$ is defined over $K_+(\n)$, absolutely irreducible, 
and isomorphic to $\Gamma(\n) \backslash \Omega^r$ over $C$ (this also follows
form the theory of modular forms, in particular, calculating the field 
of definition of the Fourier coefficients of Eisenstein series of weight one).

We prove a lemma about rank one Drinfeld modules which will be used later on. Let
$\C$ be the Carlitz module given by $\C_T = T + \tau$. 

\vspace{1ex}
\paragraph\label{lemma-Carlitz} {\bf Lemma.} \ {\sl For any rank one module $\psi$, its $\n$-torsion $\psi[\n]$ is isomorphic
as $G_K$-module to $\C[\n]$, where $\C$ is the Carlitz module.}

\vspace{1ex} 
\pf There
exists an isomorphism $v \in \bar{K}$ with $v\psi v^{-1}=\C$, and this maps a $\n$-torsion 
point $x \in \psi[\n]$ to $vt \in \C[\n]$. Now the $G_K$-action factors
through $K(\n)$, and the action of $\a \in (A/\n)^* = \Gal (K(\n)/K)$ is
given by $t \mapsto \psi_\a(t)$ and $t \mapsto \C_\a(t)$ respectively. From
this, $G_K$-equivariance follows immediately. \qed

\vspace{1ex}
\paragraph\label{T-module} {\bf $T$-motives} (Anderson, \cite{Anderson:86}).  \
Fix an $A$-Drinfeld module $\phi$. The skew $L$-algebra 
$L\{\tau \}$  becomes  a
$B=L\otimes _{\F_q} A$-algebra, denoted $M_\phi$, if we let
$$\kappa \otimes a \cdot m = \kappa \cdot m \circ \phi _a \mbox{
 and }m (\kappa \otimes a) = \kappa^q  \cdot m \circ \phi _a,$$
for every $m \in M_\phi,\; a \in A $ and $\kappa \in L.$ We note that $B$ is a polynomial 
 ring in one variable over $L:$
$B =L[Y]$ where $Y=1\otimes T;$ and $M_\phi$ is a free $B$-module of rank $r$
with  $\{1,\tau , \tau^2 ,\cdots , \tau^{r-1} \}$ as a basis.  We call $M_\phi$ the {\sl $T$-motive} associated to $\phi$. The category of $T$-motives admits tensor products,
so it makes sense to speak of the highest ($r$-fold) exterior power (=determinant)
$\ext M_\phi$ of $M_\phi$ as a $B$-module. It is a free rank one $B$-module with 
as basis $\{1 \land \tau \land \cdots \land  \tau ^{r-1}\}$; so it is again a $T$-motive, induced
by a rank one Drinfeld module which we denote by  $\land \phi$. Actually, one can calculate that $\land \phi_T = T + (-1)^{r-1} a \tau$ if $\phi$ has rank $r$ and $a$ is the leading term of $\phi_T$. (cf.\ \cite{Heiden:01}).
The determinant is 
compatible with principal level structures in the sense that $ (\wedge \phi)[\n] \cong \wedge (\phi[\n]) $
as $A[G_K]$-modules. 
It implies in particular the following

\vspace{1ex}
\paragraph\label{det-rep} {\bf Lemma.} \ {\sl Let $\rho_{\phi,\n,K} : G_K \rightarrow GL(r,A/\n)$ be the Galois representation associated to the $\n$-torsion of a rank $r$ Drinfeld module $\phi$. Then $\det \rho_{\phi,\n,K} : G_K \rightarrow (A/\n)^*$ is the Galois representation associated to 
the $\n$-torsion of the rank one Drinfeld module $\land \phi$; in particular, 
the action is given by class field theory. \qed}

\vv

\vv

\sectioning\label{X0} \ {\bf The scheme $Y_0^r(\n)$ and its twist}

\vv

We will now introduce the higher dimensional moduli space 
$Y_0^r(\n)$ and describe an Atkin-Lehner style involution that acts on it.

\vspace{1ex}
\paragraph\label{X0n} {\bf The scheme $Y_0^r(\n)$.} \ We consider the moduli functor over {\sf Sch}$/A$ associated to the moduli problem of pairs $(\phi,F)$ over an $A$-scheme $S$,
where $\phi$ is a rank $r$ Drinfeld module over $S$ and $F$ is a full flag of
$S$-subgroup schemes of the $\n$-torsion of $\phi$, i.e., 
\begin{eqnarray*}
& & F = (\phi[\n] = F_r \supset F_{r-1} \supset \dots \supset F_1 \supset F_0 = 0) \\ & & \mbox{with } F_i/F_{i-1} \times_S T \cong \n^{-1}A/A \times_S T, \ \forall i, 
\end{eqnarray*}
where $T \rightarrow S$ is a finite \'etale base change.
Isomorphism is given by an isomorphism of Drinfeld modules that induces isomorphisms on the different $F_i$ for all $i$. This functor is coarsely represented by a quasi-projective scheme $M_0^r(\n)$ of finite type over $A$, since it is actually the quotient of $M^r(\n)$ by the action of the upper triangular matrices in 
$GL(2,A/\n)$. Since $(A/\n)^*$ (as diagonal matrices) belongs to this subgroup, the components of
$M^r(\n)(C)$ are acted upon transitively, so that $M^r_0(\n)$ is absolutely irreducible. We will denote it by $Y_0^r(\n)$.

We now focus on the case $\n=T$, which is the only one that will be important for us.

\vspace{1ex}
\paragraph\label{X0nL} {\bf Lemma.} \ {\sl $Y^r_0(T)$ is an open reduced subscheme of $\P^{r-1}$, in particular, $Y^r_0(T) \times L$ for any $K$-field $L$ has a rational function field.}

\vspace{1ex}
\pf We can restrict the functor to {\sf Sch}$/K$. The scheme we have in mind in the theorem (which we call $M$) is the complement of the variety cut out by the equation $\prod a_i = 0$ in $\P^{r-1}$, where $a_i$ are coordinates on projective space. First of all, we have to check that there is a bijection between the $C$-points of $Y^r_0(T)$ and that scheme. This follows from the following considerations: if $(\phi,F)$ is a point of $Y_0^r(T)(C)$, then we can write
$$ \phi_T = (a_1 \tau +T)(a_2 \tau + 1) \dots (a_r \tau +1) $$
for $a_i \in C^*$ such that $F_i = \ker (a_{r-i+1} \tau + \dots) \dots (a_r \tau + 1)$. We map this point to $(a_1:\dots:a_r)$ in $\P^{r-1}$. The map is 
well-defined. Indeed, if $u \in C^*$ is an isomorphism $(\phi,F) \cong (\phi',F')$ (corresponding to $a'_i$), then $$ \phi'_T = u \phi_T u^{-1} = (u^{1-q} a_r \tau + T)\dots (u^{1-q} a_1 \tau + 1). $$
Then $u (F_1) = F'_1$, so that $u^{q-1} a_1 = a_1´$. Next, we have $u(F_2) = F'_2$; since $C\{\tau\}$ has a right division algorithm (Ore, cf.\ \cite{Goss:96}), this implies 
$u^{q-1} a_2 = a'_2$, etc. In the end, the vectors $(a_i)$ and $(a'_i)$ are identical up to a scalar factor $u^{q-1}$. 

Secondly, if $(a_i)$ is a point of ${\bf P}^{r-1}$ such that $\prod a_i \neq 0$, then it
clearly corresponds to a point of $Y_0^r(T)$. This is because knowing $\phi_T$ fixes the Drinfeld module (so there are no algebraic relations between the coordinates $a_i$ -- this is wrong for higher level).  This construction gives a section of the natural transformation between our functor and Mor$_M$, 
so if we are given a natural transformation $\Psi'$ between our functor and Mor$_{M'}$ for another scheme $M'$, then we can define a map from Mor$_M$ to Mor$_{M'}$ by postcomposing this section with $\Psi'$. This shows 
that $M$ represents $Y^r_0(T)$. \qed

\vspace{1ex}
\paragraph\label{X0nC} {\bf Remark.}  \ Picking a basis for $\phi[T]$ over $C$ allows one to represent the moduli space $Y^r_0(\n) \times C$ as the quotient of
Drinfeld´s $(r-1)$-dimensional symmetric space $\Omega^r$ by the arithmetic 
subgroup of $GL(r,A)$ given by matrices congruent to an upper triangular matrix
modulo $\n$, but we will not need this description. 

\vv

\vspace{1ex}
\paragraph\label{definvo} {\bf The involution $w$.} \ Assume that $r=2s$ is even. We define the Atkin-Lehner involution $w=w_T$ on $Y_0^r(T)$ to be the map that maps the point $$(\phi_T=(a_1 \tau +T)(a_2 \tau + 1) \dots (a_r \tau +1),F)$$ to its image under the isogeny $$u=(a_{s+1}\tau+1) \dots (a_r\tau+1).$$ One computes immediately that
\begin{eqnarray*}
u \phi_T u^{-1} &=& (a_{s+1} \tau + 1) \dots (a_r \tau + 1)(a_1 \tau + T) \dots (a_s \tau + 1) \\
&=& (a_{s+1} \tau + 1) \cdot  T \cdot T^{-1} \dots (a_r \tau + 1)\cdot T \cdot T^{-1} (a_1 \tau + T) \dots (a_s \tau + 1) \\
&=& (a_{s+1} T^q \tau + T) \dots (a_r T^{q-1} \tau + 1)(T^{-1} a_1 \tau + 1) \dots (a_s \tau + 1),
\end{eqnarray*} 
so that 

\vspace{1ex}
\paragraph\label{leminvo} {\bf Lemma.} {\sl With respect to the natural embedding $Y_0^r(T) \rightarrow {\bf P}^{r-1} : (\phi,F) \rightarrow (a_i)$ of the previous lemma, the involution $w$ maps 
$(a_i)$ to $$ (T^qa_{s+1}:T^{q-1} a_{s+2}:T^{q-1} a_{s+3}:\dots:T^{q-1} a_r:T^{-1} a_1:a_2:a_3:\dots:a_s). \ \Box$$}
In this way, one sees directly that $w$ is an involution and that it induces an
automorphism of $Y_0^r(T)$. 

\vspace{1ex}
\paragraph\label{defYY} {\bf The twisted moduli space $\YY$.} \  We now assume that $L/K$ is a quadratic field with non-trivial $K$-automorphism $\sigma$, and we look at the twisted moduli space
$$ \YY := (w \times \sigma) \backslash (Y_0^r(T) \times_{\Spec\ K} \Spec \ L), $$
i.e., the quotient on which $w$ acts like $\sigma$.

\vspace{1ex}
\paragraph\label{YYrat} {\bf Lemma.} {\sl $\YY$ is a rational variety
over $K$ of dimension $r-1$.}

\vspace{1ex}
\pf The action of $w$ extends to $\P^{r-1}$ in the obvious way described by the previous lemma. We then set 
$$ P := (w \times \sigma) \backslash (\P^{r-1}_K \times_{\Spec\ K} \Spec \ L). $$
Since $Y_0^r(T) \times K$ is dense open in $\P^{r-1}_K$, $\YY \times K$ is dense open in $P$, 
an in particular, they have the same function field. 

$P$ is by definition a Brauer-Severi variety (cf.\ \cite{Jahnel:00}) with splitting field $L$. Indeed, over $L$, it is isomorphic to $w \backslash \P^{r-1}_L$, which in its turn is isomorphic to $\P^{r-1}_L$ (since $\langle w \rangle$ is an abelian group of order coprime to $p$ acting linearly on projective space, one can assume that it acts by
characters (of order $\leq 2$) on suitable variables of the homogeneous coordinate ring $R=L[x_1,\dots,x_r]$ of $\P^{r-1}_L$, so that $R^{\langle w \rangle} = L[x_1,\dots,x_l,x_{l+1}^2,\dots,x_r^2]$ for $l$ the dimension of the invariant subspace -- which, by the way, equals $s$).
 
Hence to prove that $P$ is $K$-isomorphic to $\P^{r-1}$, it suffices to 
find a point in $P(K)$ (Ch\^{a}telet, cf.\ \cite{Jahnel:00}, 4.8). For example, the point  $$(T^{\frac{q+1}{2}} \cdot a_{s+1}:0:\dots :0:a_{s+1}:0:\dots :0)$$ is a $K$-rational fixed point of 
$w$ on $\P^{r-1}$ and hence descends to a $K$-rational point of $P$.
This proves that $\YY$ also has function field $K(t_1,\dots,t_{r-1})$. \qed 

\vv

\vv

\sectioning \label{Galrep} {\bf A Galois representation associated to rational points of $\YY$.}

\vv

\paragraph\label{Krat} {\bf Rational points of $\YY$.} \ A $K$-rational point of $\YY$ is by definition a pair $(\phi,F)$ defined over $L$ such that
$$(\phi^\sigma,F^\sigma) = (u\phi u^{-1},u(F)).$$ In particular, $\phi$ is a
Drinfeld module over $L$ which admits a $T$-isogeny (i.e., with kernel
a subspace of $\phi[T](L)$) of rank $s$ to 
its Galois conjugate $\phi^\sigma$. Note furthermore that $w (\phi[T])= (\ker\ u)^\sigma$, since both are equal to $\ker (\hat{u})$ for $\hat{u}$ the isogeny
dual to $u$ (i.e., $\hat{u} u = \phi_T$). 

\vspace{1ex}
\paragraph\label{fieldL} {\bf Choice of the field $L$.} \ Fix a prime $\fp$ of degree $d>1$ in $A$. From now on, we set $L=K(\sqrt{\fp})$. We claim that $L$ is a subfield of $K(\fp)$. Indeed, $K(\fp)$ is the field generated by the $\fp$-torsion 
of the Carlitz module, i.e., the roots of $\C_\fp$ if $\C_T=TX+X^q$. We calculate
the discriminant $D$ of $\C_{\fp}$ as the resultant of $\C_\fp$ and its $X$-derivative $\C´_\fp=\fp$, to find $D=\fp^{q^d}$. Hence $L=K(\sqrt{D})$, which clearly is a subfield
of $K(\fp)$. 

\vspace{1ex}
\paragraph\label{rho} {\bf The associated Galois representation.} \ Given such a pair $(\phi,F)/L$ corresponding to a rational point of $\YY$, we construct a Galois representation 
$$ \rho : G_K \rightarrow PGL(r,A/\fp) $$
as follows. Recall that $G_L$ has index 2 in $G_K$, the quotient being generated by $\sigma$ (the generator of of $\Gal(L/K)$). If $\tau \in G_L$, then we set
$ \rho(\tau) = \rho_{\phi,\fp,L}(\tau)$ modulo scalars in $GL(r,A/\fp)$. On the other hand, the isogeny $u$ induces an isomorphism $\phi[\fp] \cong \phi[\fp]^\sigma$ since $T$ and $\fp$ are coprime (the dual does the same in the other direction), so it makes sense to set
$ \rho(\sigma)(x) = \hat{u}(x^\sigma)$ as an automorphism of $\phi[\fp]$, i.e.,
an element of $GL(r,A/\fp)$, which we then again consider modulo scalars.

\vspace{1ex}
\paragraph\label{welldef} {\bf Lemma.} \ {\sl If $\mbox{{\rm Aut}}(\phi)={\bf F}_q^*$, this $\rho$ is a well-defined
group representation of $G_K$.}

\vspace{1ex}
\pf We have to check that $\rho(\tau_1 \tau_2) = \rho(\tau_1) \rho(\tau_2)$ for all $\tau_1, \tau_2$ in $G_K$. It is clear if both $\tau_i$ belong to $G_L$. 

Now assume $\tau_1 \notin G_L, \tau_2 \in G_L$, so $\tau_1 \tau_2 \notin G_L$.
Then if $x \in \phi[\fp]$, $x^{\rho(\tau_1 \tau_2)}$ is by definition 
$\hat{u}(x^{\tau_1 \tau_2})$, whereas  $$x^{\rho(\tau_1) \rho(\tau_2)} = (\hat{u}(x^{\tau_1})^{\tau_2}) = \hat{u}(x^{\tau_1 \tau_2}).$$ 
The case $\tau_1 \in G_L, \tau_2 \notin G_L$ is similar. 

Finally, if $\tau_1 \notin G_L, \tau_2 \notin G_L$, then $\tau_1 \tau_2 \in G_L$, so on the one hand,  $x^{\rho(\tau_1 \tau_2)} = x^{\tau_1 \tau_2}$, and on the other hand,  $$x^{\rho(\tau_1) \rho (\tau_2)}
= \hat{u}(\hat{u}(x^{\tau_1}))^{\tau_2} = \hat{u} \hat{u}^\sigma (x^{\tau_1 \tau_2}).$$ We claim that $\hat{u}^\sigma = \lambda u$ for some $\lambda \in \F_q^*$.
Taking this claim for granted, we conclude (using $\hat{u} u = \phi_T$) that $$x^{\rho(\tau_1) \rho (\tau_2)}= \lambda T x^{\tau_1 \tau_2}$$ where the multiplication with 
$\lambda T$ is that of a scalar matrix modulo $\fp$ via the canonical identification $\phi[\fp]=(A/\fp)^r$; hence after dividing out scalars, we indeed get $x^{\tau_1 \tau_2}$.

Now for the proof of the claim, we observe that $u^\sigma$ and $\hat{u}$ are
two $T$-isogenies from $\phi^\sigma$ to $\phi$ defined over $L$.  We have remarked in \ref{Krat} that the two isogenies have the same kernel, so using the right division algorithm in $L\{\tau\}$, we see that they differ by a scalar 
element in $ \lambda \in L^*$. On the other hand, we compute (using that $\sigma$ has order two and that $u$ is an isogeny) for any $a \in A$:
$$ \lambda \phi_a^\sigma u = \lambda u \phi_a  = \hat{u}^\sigma \phi_a = (\hat{u} \phi_a^\sigma)^\sigma = (\phi_a \hat{u})^\sigma = \phi_a^\sigma \hat{u}^\sigma = \phi_a^\sigma \lambda u.$$
Canceling $u$, this says that $\lambda$ is an endomorphism of $\phi^\sigma$, which 
together with the fact that it is a scalar implies $\lambda \in \mbox{Aut}(\phi^\sigma)$. 
Finally, our assumption that Aut$(\phi)={\bf F}_q^*$ gives the claim. \qed

As a preparation for the statement of the next lemma, note that if $(r;q^d-1)=2$, then $(A/\fp)^*/(A/\fp)^{*r} = (A/\fp)^*/(A/\fp)^{*2} \cong {\bf Z}/2 \cong  \{ \pm 1 \}$.

\vspace{1ex}
\paragraph\label{detimage} {\bf Lemma.} \ {\sl Assume that Aut$(\phi)={\bf F}_q^*$, that $r$ is even and $(r;q^d-1)=2$. Then $\det(\rho(L))=1$ as an element of $(A/\fp)^*/(A/\fp)^{*r} \cong \{ \pm 1 \}$, and $$\det \rho (\sigma)=-[ \frac{\zeta T}{\fp}],$$ where $[\frac{\cdot}{\fp}]$ denotes the quadratic character modulo $\fp$ and $\zeta \in \F_q^*$ is a constant depending only on $r$.}

\vspace{1ex}
\pf We compute the action of Galois on the determinants of the occurring Drinfeld modules ($\land \phi$ and $\land \phi^\sigma$) via \ref{det-rep}. For an element $\tau \in G_L$, $\det \rho (\tau)$ (acting on the Carlitz module, which we can assume by lemma \ref{lemma-Carlitz}) equals $\det \rho_{\phi,\n,L}(\tau)$, which is a square
in Gal$(K(\fp)/K) = (A/\fp)^*$ since $L$ is the subfield of $K(\fp)$ of degree two over $K$. So we only have to compute $\det \rho (\sigma)$.

Using lemma \ref{lemma-Carlitz}, $\sigma: \land \phi[\fp] \rightarrow \land \phi^\sigma[\fp]$ is equivalent to the action of $\sigma$ on $\C[\fp]$, which is given by class field theory as an element of order two in $(A/\fp)^*$ (since $L$ lies as a quadratic
subfield in $K(\fp)$).

For the action of $\hat{u} : \land \phi^\sigma \rightarrow \land \phi$, we
rather look at the induced ($B$-algebra) homomorphism between associated 
$T$-motives $\hat{u} : M_{\land \phi^\sigma} \rightarrow M_{\land \phi}$.
On both sides, we pick the basis $e=1 \land \dots \tau^{r-1}$. Since
the modules are one-dimensional, we know that $\hat{u}$ acts like 
multiplication with an element $1 \otimes a$ of $B$. To determine this element,
we observe that $u \hat{u}(=\phi_T)$ acts like $1 \otimes T$ on $M_\phi$; hence
it acts like $1 \otimes T^r$ on $\land \phi$. Clearly, $\deg_T a = \deg_\tau u = s$, so $a=\zeta T^s$ for some $\zeta \in \F_q^*$. 

Taking these two calculations together, we see that $\det \rho$ acts
like multiplication with the product of a non-square (coming from $\sigma$) and $\zeta T^s$. Since
$s$ is odd we get that modulo squares, the determinant equals $-[\frac{\zeta T}{\fp}]$.

From this description, it is not yet clear that $\zeta$ does not depend on 
the particular ``point'' $(u: \phi \rightarrow \phi^\sigma)$ of $\YY$. This, however, can be seen as follows. Suppose that we calculate formally, starting with a Drinfeld module $$\phi_T=(a_1 \tau +T)(a_2 \tau + 1) \dots (a_r \tau +1)$$ over the function field $L'=L(a_i)$ of $\YY$. The isogeny $u$ is (on the level of motives) a $B$-module homomorphism between $M_\phi$ and $M_{\phi^\sigma}$, and if we express $u$ as a matrix w.r.t.\ the basis $\{1,\dots,\tau^{r-1}\}$ on both sides, it has entries in $B':=L' \otimes A$. Hence the same holds for its determinant, which we already know is of the form $l \otimes T^s$ for some $l \in L'$. We have expressed $\det(u)$ w.r.t.\ the bases $1 \wedge \dots \wedge \tau^{r-1}$ of $M_{\wedge \phi}$ and $M_{\wedge \phi^\sigma}$. Now we can make both of these one-dimensional $B'$-modules into the Carlitz module by an isomorphism, which here is of the form $i \otimes 1 \in B'$ as a $B'$-module homomorphism. Now we see that $\det(u)$ w.r.t.\ these bases is given by  $li \otimes T^s$, which we know from before equals $\zeta \otimes T^s$. Identifying both results, we see that $\zeta$ depends as a rational
function on the coordinates $a_i$ of the variety $\YY$. This means that the
{\sl finite} stratification of $\YY$ induced by what $\zeta$ occurs in the corresponding Galois representation is actually by {\sl subvarieties}. Since there are only finitely many, at least one of them is equidimensional with $\YY$, but
$\YY$ is absolutely irreducible, so it equals this subvariety. It follows that $\zeta$ is constant on 
all of $\YY$ (namely, what it is on that subvariety). \qed

\vv

\vspace{1ex}
\paragraph\label{exdet} {\bf Example.} \ One can make the considerations in the lemma very explicit if one fixes $r$. For example, if $r=2$, write $$\phi_T = (a_1 \tau + T)(a_2 \tau + 1)=a_1 a_2^q \tau^2 + (Ta_2+a_1) \tau + T$$ and $u=a_2 \tau + 1$. W.r.t.\ the basis $\{1,\tau\}$, $u$ is the matrix
$$ u = \left( \begin{array}{cc} 1 &  a_1^{-1}(1 \otimes T - T) \\ a_2 & -Ta_1^{-1}a_2 \end{array} \right), $$
so the determinant is $ \det(u) = -a_1^{-1} a_2 \otimes T.$ One calculates
that $$\wedge \phi_T = T - a_1 a_2^q \tau \mbox{ and }\wedge \phi_T^\sigma = \wedge u \phi_T u^{-1} = T - a_2 a_1^q \tau,$$ so that \begin{eqnarray*} & & \mbox{the isomorphism } {\mathcal C} \cong \wedge \phi_T \mbox{ is multiplication by }(a_2 \sqrt[q-1]{-a_1 a_2})^{-1} \otimes 1 \\  & & \mbox{the isomorphism } \wedge \phi^\sigma_T \cong {\mathcal C} \mbox{ is multiplication by } a_1 \sqrt[q-1]{-a_1 a_2} \otimes 1. \end{eqnarray*} Finally, we get that as an automorphism of the Carlitz module, $\det(u)=-1 \otimes T$ (so here, $\zeta$ is not a square in ${\bf F}_q^*$). 

\vv

\vv

\sectioning \label{Hilb} {\bf Proof of the main theorem.}

\vv

The moduli functor, and hence also the moduli space $M^r(\fp)$, comes with a natural action of $GL(r,A/\fp)$ and the quotient of this action is $M^r(1)$. So $\rho_{{\phi},\fp,K(M^r(1))}$ is surjective for $F$ the function field of $M^r(1)$ and
$\phi$ the ``universal'' Drinfeld module over $M^r(1)$. Since $\YY \rightarrow M^r(1)$ is finite, the same holds if
we replace $K(M^r(1))$ by $K(\YY)$.

If we replace $K$ by $K(\YY)$ and $L$ by $L(\YY)$, the construction of \S \ref{Galrep} provides us with a representation $$\rho \ : \ G_{K(\YY)} \rightarrow PGL(r,q^d)$$ whose image we will now determine. The arguments in lemma \ref{detimage} show that $\det \rho = 1$
if the generic $\phi$ does
not have non-trivial automorphisms and if we can choose $\fp$ such that $[\frac{\zeta T}{\fp}]=-1$. That the first is true is obvious, since Aut$(T+a_1 \tau + \dots + a_r \tau^r) \neq {\bf F}_q^*$ implies the non-generic condition that
$a_1 = 0$. That we can accomplish the second condition comes out of the following
lemma:

\vspace{1ex}
\paragraph\label{HM} {\bf Lemma.} \ {\sl For any $d>1$, there exists 
an irreducible polynomial $\fp$ of degree $d$ in $A$ such that $[\frac{T}{\fp}]$ assumes a given value ($\pm 1$).}

\vspace{1ex}
\pf By quadratic reciprocity in function fields, it suffices to show that $[\frac{\fp}{T}]$ can attain a given value. But this depends on the fact whether or not the constant term of $\fp$ is a square in $\F_q$. Hence it 
certainly suffices to show that for any given $\xi \in \F_q^*$, there exists
an irreducible polynomial $\fp$ of degree $d$ over $K$ whose constant term equals $\xi$; for a proof of this, see Hansen-Mullen (\cite{Hansen:92}, p.\ 642). \qed

\vv

Since $\rho_{{\phi},\fp,K(\YY)}$ is surjective, we get that 
$$ \rho \ : \ \ G_{K(\YY)} \rightarrow PSL(r,q^d)$$
is surjective, and by lemma \ref{YYrat}, we find that 
$K(\YY) = K(t_1,\dots,t_{r-1}).$ 

Finally, the fixed field of the separable
closure of $K(\YY)$ under the kernel of $\rho$ gives a $PSL(r,q^d)$-extension of $K(\YY)$, 
which is regular, since the non-regular subextension $K_+(\fp)$ of $K(M^r(\fp))$ corresponds to the (inverse image of) the center of $GL(r,q^d)$, which 
intersects $\ker(\rho)$ trivially (since the determinant of $\rho$ is 
trivial). This finishes the proof of the theorem. \qed
 
\vv

\vv

\sectioning \label{explicit} {\bf Proof of the proposition --- an explicit example.}

\vv

\vspace{1ex}
\paragraph\label{hauptmoduln} {\bf Setup.} \ Let $N$ be irreducible in $A$ of degree 2. We will use a method similar to the one in \cite{Shih:78}, \S 5 to make the construction explicit --- in more complicated cases the method of \cite{Malle} should be applied. We will assume here that the reader is familiar with the basic notions from the theory of Drinfeld modular forms, as explained for example in \cite{Gekeler:88}.  We make the convention that all of our modular forms and functions are normalized so the leading term of their expansion at the infinite cusp is $=1$. The curves $X_0(T)$ and $X_0(N)$ are rational, and good {\sl Hauptmoduln} for them are given by 
$$ f_T(z) := {\frac{h(Tz)}{h(z)}} \mbox{ and }  f_N(z) := \sqrt[q+1]{\frac{h(Nz)}{h(z)}} $$
respectively, where $h$ is a Poincar\'e series of weight $q+1$ for $GL(2,A)$. The curve $X_0(NT)$ is of genus $q$, but its quotient $X_+(NT)$ by the Atkin-Lehner involution $w_{NT}$ is rational. Note that for $a \in A$, $w_a$ acts like $w_a(z)=-\frac{1}{az}$ on the upper half plane. Also note that for the uniformizer $s(z):=t(z)^{q-1}$ of $X(1)$ at the infinite cusp, we have $$t(az)=\frac{t^{|a|}}{{\mathcal C}_a(t(z)^{-1})t^{|a|}},$$ where $\mathcal C$ is the Carlitz module. In particular, the order of $f(az)$ at $\infty$ is $|a|$ times the order of $f$ at $\infty$.

\vspace{1ex}
\paragraph\label{hauptmoduln+} {\bf Lemma.} {\sl A {\rm Hauptmodul} for $X_+(NT)$ is given by  
$$ f(z) := \sqrt[(q^2-1)]{\frac{h(NTz)h(z)}{h(Tz)h(Nz)}}. $$}

\pf Observe that the quantity underneath the root sign is clearly $w_{NT}$-invariant, and has divisor on $X_+(NT)$ supported at the cusps, of which there are exactly two. Hence we can extract a $d$-th root of this function in the function field of $X_+(NT)$ for $d$ precisely equal to its order at, say, the cusp at infinity. Now since $h=t(1+O(s))$, the expression under the root is $s^{q^2-1}+$ higher terms, whence the result. \qed

\vspace{1ex}
\paragraph\label{fwt} {\bf Lemma.} {\sl $w_T f = f^{-1}.$} 

\vspace{1ex}
\pf The divisor of $w_T f$ is the image of the divisor of $f$ under $w_T$. As $w_T$ interchanges the two cusps of $X_+(NT)$, we get that $\mbox{div}(w_Tf)=-\mbox{div}(f)$. Hence $w_T f$ differs from $f^{-1}$ by a constant, which is seen to be $=1$ by expanding.  \qed 

\vspace{1ex}
\paragraph\label{fields} {\bf Function fields of $\tilde{Y}$'s.} \ The function field of $X_0(NT)$ is generated by $f$ and $f_T$. As $w_T f_T = \frac{1}{T^{q+1}f_T}$, the function field of $\tilde{Y}(T)$ is generated by $x$ and that of 
$\tilde{Y}(NT)$ by $x$ and $y$, where 
$$ x = \sqrt{N} \frac{T^{\frac{q+1}{2}}f_T-1}{T^{\frac{q+1}{2}}f_T+1} \mbox{ and } y = \sqrt{N} \frac{f-1}{f+1}. $$
Note that $X(\mbox{sth})$ is the Galois closure of $X_0(\mbox{sth})$, so our theoretical result implies that for suitable $N$, the 
splitting field of the minimal polynomial of $y$ over $K(x)$ is a $K$-regular Galois extension 
with Galois group $PSL(2,q^2)$. Since $\zeta$ is a non-square for $r=2$ (cf.\ \ref{exdet}), we should choose $[\frac{T}{N}]=1$, i.e., the constant term of $N$ should be a square modulo $T$ (by reciprocity). 

To find an explicit equation, it suffices to find a relation between $f$ and $f_T$, and solve it for $x$ and $y$. Now one calculates that 
$$ w_{NT} f_T =\frac{1}{T^{q+1} f_T f^{q^2-1}}, $$
so 
$$ T^{q+1} f_T f^{q^2} + \frac{f}{f_T} $$
 is invariant under $w_{NT}$ and hence belongs to $K(f)$. We see in particular by looking at the order of poles and zeros that if is in fact a polynomial in $f$ of degree $q^2+1$. This proves the proposition. \qed 

\vspace{1ex}
\paragraph\label{expl-eqs} {\bf Explicit equations.} To find this polynomial, it suffices to compare series expansions, and for this one can use the fact
that 
$$ h=t(-\frac{1}{U_1}+s^{q^3-q^2}(\frac{1}{U_1}-\frac{1}{U_1^2})-s^{q^3-1}+O(s^{q^3})), $$
where $U_1=1-s^{q-1}+(T^q-T)s^q$ (cf.\ Gekeler, \cite{Gekeler:88}, (10.4)). 
The example from the introduction was computed setting $q=3, N=T^2+1$ in {\tt gp-pari} and {\tt Maple}. 

\vv

\vv

{\footnotesize

\noindent {\bf Acknowledgments.} The first author is honorary fellow of the Fund for Scientific Research - Flanders (FWO -Vlaanderen). The second author is partially supported by IKY (Greece). Part of this work was done at the MPIM. We thank Gert-Jan van der Heiden for useful comments.  

\vv

\vv

\bibliographystyle{amsplain}

\begin{thebibliography}{10}

\bibitem{Abhyankar:95}
S.S.~Abhyankar, \emph{Mathieu group coverings and linear group coverings}, in: Recent developments in the inverse Galois problem (Seattle, WA, 1993), 293--319, Contemp.
Math., 186, Amer. Math. Soc., Providence, RI, 1995.

\bibitem{Abhyankar:00b}
S.S.~Abhyankar and P.H.~Keskar, \emph{Descent principle in modular Galois theory}, Proc. Indian Acad. Sci. Math. Sci. \textbf{111} (2001), no. 2,
139--149. 

\bibitem{Anderson:86}
G.W.~Anderson, \emph{t-Motives}, Duke Math.\ J.\ {\bf 53} (1986), 457--502.

\bibitem{Gekeler:88} 
E.-U.\ Gekeler, \emph{On the coefficients of Drinfeld modular forms}, Invent. Math. \textbf{93} (1988), no. 3, 667--700.

\bibitem{Goss:96}
D.~Goss, \emph{ Basic structures of function field arithmetic}, 
Ergebn.\ Math.\ Grenzg.\ (3) vol.\ 35, 
Springer-Verlag, Berlin, 1996.

\bibitem{Hansen:92}
T.~Hansen and G.L. Mullen, \emph{Primitive polynomials over finite fields}, Math.\ Comp.\ {\bf 59} (1992), 639--643.

\bibitem{Heiden:01}
G.-J.\ van der Heiden, \emph{Weil-pairing for Drinfeld modules of rank 2},
preprint, Universiteit Groningen, 2001.

\bibitem{Jahnel:00}
J.~Jahnel, \emph{The Brauer-Severi variety associated with a central simple
algebra: a survey}, Lin.\ Alg.\ Grps.\ and Related Structures (electronic) \textbf{52} (2000), 1--60.

\bibitem{Malle}
G.~Malle, \emph{Polynome mit Galoisgruppen ${\rm PGL}\sb 2(p)$ und ${\rm PSL}\sb 2(p)$ \"uber $\Q(t)$}, Comm. Algebra \textbf{21} (1993), no. 2, 511--526. 

\bibitem{Nori}
M.V. Nori, \emph{Unramified coverings of the affine line in positive characteristic}, in: Algebraic Geometry and its applications. Collection of
papers from S.\ Abhyankar's birthday conference (West Lafayette, 1990) (C.L. Bajaj, ed.), 209-212, Springer-Verlag, 1994. 

\bibitem{Pop}
F.~Pop, \emph{Embedding problems over large fields}, 
Ann. of Math. (2) \textbf{144} (1996), no. 1, 1--34. 

\bibitem{Shih:78}
K.-y. Shih, \emph{$p$-division points on certain elliptic curves}, Compos.\
Math.\ {\bf 36} (1978), 113--129. 

\bibitem{Vladuts:91}
S.~Vladut, \emph{Kronecker's Jugendtraum and modular functions}, Studies in 
the Development of Modern Mathematics, Vol.\ 2, Gordon and Breach, New York, 1991.

\end{thebibliography}
\providecommand{\bysame}{\leavevmode\hbox to3em{\hrulefill}\thinspace}

\vv

\vv

\noindent Department of Mathematics,
Utrecht University, P.O. Box 80010, 3508 TA  Utrecht, The Netherlands (gc)\footnote{All correspondence should be sent to this author at
this address}

\vv

\noindent Askoutsi 20, TK 74100 Crete, Greece (mt)

\vv

\noindent e-mail: {\tt cornelissen@math.uu.nl}, {\tt tripolit@csd.uoc.gr}
}

\end{document}